\documentclass[11pt]{article}
\usepackage{latexsym}
\usepackage{amssymb,amsmath,amsfonts}
\usepackage[mathscr]{eucal}
\usepackage{amscd}
\usepackage{amsfonts,cite}
\usepackage{epsfig, amsmath, amsthm, amssymb}
\usepackage{latexsym}
\usepackage{graphics,epstopdf}
\begin{document}
\thispagestyle{empty}
\parindent=0mm
\begin{center}
{\Large{\bf Symmetry identities for 2-variable Apostol type \\and related polynomials}}\\~~

{\bf Subuhi Khan${}^{a,\star}${\footnote{${}^\star$Corresponding author; E-mail:~subuhi2006@gmail.com (Subuhi Khan)}},
Mumtaz Riyasat${}^b$ {\footnote{E-mail:~mumtazrst@gmail.com (Mumtaz Riyasat)}}
and Ghazala Yasmin${}^c$ {\footnote{Email: ghazala30@gmail.com (Ghazala Yasmin)\\}}}\\

${}^{a,b}$Department of Mathematics, Aligarh Muslim University, Aligarh, India\\
${}^c$Department of Applied Mathematics, Faculty of Engineering\\
Aligarh Muslim University\\
Aligarh, India\\
\end{center}

\parindent=0mm
{\bf Abstract:}~In this article, certain symmetry identities for the 2-variable Apostol type polynomials are derived. By taking suitable
values of parameters and indices, the symmetry identities for the special  cases of the 2-variable Apostol type polynomials
are established. Further, the symmetry identities for certain members belonging to the 2-variable Apostol type polynomials are also considered.\\


\noindent
{\large{\bf 1.~Introduction and preliminaries}}\\
\vspace{.35cm}

\parindent=8mm
The class of Appell sequences is an important class of polynomial sequences and very often found in different applications in pure and applied mathematics. Properties of Appell sequences are naturally handled within the framework of modern classical umbral calculus by Roman \cite{Roman}. The classical Bernoulli polynomials $B_n(x)$, the classical Euler polynomials $E_n(x)$ and the classical Genocchi polynomials $G_n(x)$, together with their familiar generalizations $B_n^{(m)}(x)$, $E_n^{(m)}(x)$ and $G_n^{(m)}(x)$ of (real or complex) order $m$, belong to the class of Appell sequences.\\

The polynomials $B_n^{(m)}(x)$, $E_n^{(m)}(x)$ and $G_n^{(m)}(x)$ are defined by the following generating functions \cite{Erd3,Sandor}:
$$\Big(\frac{t}{ e^{t}-1}\Big)^m e^{xt}=\sum\limits_{n=0}^{\infty} B_n^{(m)}(x)\frac{t^{n}}{n!},~~|t| < 2\pi,\eqno(1.1)$$
$$\Big(\frac{2}{ e^{t}+1}\Big)^m e^{xt}=\sum\limits_{n=0}^{\infty} E_n^{(m)}(x)\frac{t^{n}}{n!},~~|t| < \pi,\eqno(1.2)$$
$$\Big(\frac{2t}{ e^{t}+1}\Big)^m e^{xt}=\sum\limits_{n=0}^{\infty} G_n^{(m)}(x)\frac{t^{n}}{n!},~~|t| < \pi.\eqno(1.3)$$

It is easy to see that $B_n(x)$, $E_n(x)$ and $G_n(x)$ are given, respectively, by
$$B_n^{(1)}(x)=B_n(x);~E_n^{(1)}(x)=E_n(x);~G_n^{(1)}(x)=G_n(x),~n \in \mathbb{N}_0:=\mathbb{N}\cup \{0\}.\eqno(1.4)$$

Some interesting analogues of the classical Bernoulli and Euler polynomials were first investigated by Apostol \cite{Apos} and further studied by Srivastava \cite{Srivastava}.\\

Luo and Srivastava \cite{Luo1} introduced the Apostol Bernoulli polynomials of order $m \in \mathbb{N}_0$, denoted by $\mathfrak{B}_n^{(m)}(x;\lambda),~\lambda \in \mathbb{C}$, which are defined by the generating function
$$\Big(\frac{t}{\lambda e^{t}-1}\Big)^m e^{xt}=\sum\limits_{n=0}^{\infty} \mathfrak{B}_n^{(m)}(x;\lambda)\frac{t^{n}}{n!},~~|t| < 2\pi, ~\textrm{when}~\lambda=1;~|t|<|\textrm{log} \lambda|,~\textrm{when}~\lambda\neq1,\eqno(1.5)$$
with
$$\mathfrak{B}_n^{(m)}(x;1)=B_n^{(m)}(x).\eqno(1.6)$$

The Apostol Euler polynomials of order $m \in \mathbb{N}_0$, denoted by $\mathfrak{E}_n^{(m)}(x;\lambda),~\lambda \in \mathbb{C}$ are introduced by Luo \cite{Luo} and are defined by the generating function
$$\Big(\frac{2}{\lambda e^{t}+1}\Big)^m e^{xt}=\sum\limits_{n=0}^{\infty} \mathfrak{E}_n^{(m)}(x;\lambda)\frac{t^{n}}{n!},~|t|<|\textrm{log}(-\lambda)|,\eqno(1.7)$$
with
$$\mathfrak{E}_n^{(m)}(x;1)=E_n^{(\alpha)}(x).\eqno(1.8)$$

Further, Luo \cite{Luo2} introduced the Apostol Genocchi polynomials of order $m \in \mathbb{N}_0$, denoted by $\mathcal{G}_n^{(m)}(x;\lambda),~\lambda \in \mathbb{C}$, which are defined by the generating function
$$\Big(\frac{2t}{\lambda e^{t}+1}\Big)^m e^{xt}=\sum\limits_{n=0}^{\infty} \mathcal{G}_n^{(m)}(x;\lambda)\frac{t^{n}}{n!},~|t|<|\textrm{log}(-\lambda)|,\eqno(1.9)$$
with
$$\mathcal{G}_n^{(m)}(x;1)=G_n^{(m)}(x).\eqno(1.10)$$

Recently, Luo and Srivastava \cite{Luo3} introduced a unification (and generalization) of the above mentioned three families of the Apostol-type polynomials.\\

The Apostol type polynomials $\mathcal{F}_n^{(\alpha)}(x;\lambda;\mu,\nu)$ $(\alpha \in \mathbb{N}, \lambda,\mu,\nu \in \mathbb{C} )$ of order $m$, are  defined by means of the generating function
$$\Big(\frac{2^\mu~ t^\nu}{\lambda e^{t}+1}\Big)^m e^{xt}=\sum\limits_{n=0}^{\infty} \mathcal{F}_n^{(m)}(x;\lambda;\mu,\nu)\frac{t^{n}}{n!},~|t|<|\textrm{log}(-\lambda)|.\eqno(1.11)$$

In fact, from equations (1.5), (1.7), (1.9) and (1.11), we have
$$(-1)^m \mathcal{F}_n^{(m)}(x;-\lambda;0,1)=\mathfrak{B}_n^{(m)}(x;\lambda),\eqno(1.12)$$
$$ \mathcal{F}_n^{(m)}(x;\lambda;1,0)=\mathfrak{E}_n^{(m)}(x;\lambda),\eqno(1.13)$$
$$ \mathcal{F}_n^{(m)}(x;\lambda;1,1)=\mathcal{G}_n^{(m)}(x;\lambda).\eqno(1.14)$$

Recently, Khan {\em at el} introduced the 2-variable Apostol type polynomials, denoted by $_p\mathcal{F}_n^{(m)}(x,y;\lambda;\mu,\nu)$ by combining the 2-variable general polynomials $p_n(x,y)$ \cite{KhG} as base with Apostol type polynomials $\mathcal{F}_n^{(m)}(x;\lambda;\mu,\nu)$, which are defined by the generating function of the form \cite{ATP}:

$$\Big(\frac{2^\mu~ t^\nu}{\lambda e^{t}+1}\Big)^m e^{xt}\phi(y,t)=\sum\limits_{n=0}^{\infty} {}_p\mathcal{F}_n^{(\alpha)}(x,y;\lambda;\mu,\nu)\frac{t^{n}}{n!}.\eqno(1.15)$$

Symmetry identities for the generalized Apostol Bernoulli polynomials $\mathfrak{B}_n^{(m)}(x;\lambda)$ and generalized Apostol Euler polynomials $\mathfrak{E}_n^{(m)}(x;\lambda)$ involving generalized sum of integer powers $\mathcal{S}_k(n, \lambda)$ and generalized sum of alternative integer powers $\mathcal{M}_k(n, \lambda)$ are derived in \cite{Zhang}.\\

Now, we recall the following definitions which will be used in proving the results:\\

\noindent
{\bf Definition~1.1.}~For each $k \in \mathbb{N}_0$, the sum $S_k(n)=\sum\limits_{i=0}^n i^k$ is known as the power sum defined by the generating relation
$$\sum\limits_{k=0}^{\infty}S_{k}(n)\frac{t^k}{k!}=\frac{ e^{(n+1)t}-1}{ e^{t}-1}.\eqno(1.16)$$

\noindent
{\bf Definition~1.2.}~For an arbitrary real or complex $\lambda$, the generalized sum of integer powers $\mathcal{S}_k(n, \lambda)$ is defined by the generating relation \cite{Zhang}:

$$\sum\limits_{k=0}^{\infty} \mathcal{S}_{k}(n,\lambda)\frac{t^k}{k!}=\frac{\lambda e^{(n+1)t}-1}{\lambda e^{t}-1}.\eqno(1.17)$$

We note that
$$\mathcal{S}_{k}(n,1)=S_{k}(n).\eqno(1.18)$$

\noindent
{\bf Definition~1.3.}~For each $k \in \mathbb{N}_0$, the sum $M_k(n)=\sum\limits_{i=0}^n (-1)^ki^k$ is known as the sum of alternative integer powers defined by the generating relation
$$\sum\limits_{k=0}^{\infty}M_{k}(n)\frac{t^k}{k!}=\frac{ 1-(-e^t)^{(n+1)}}{ e^{t}+1}.\eqno(1.19)$$

\noindent
{\bf Definition~1.4.}~For an arbitrary real or complex $\lambda$, the generalized sum of alternative integer powers $\mathcal{M}_k(n, \lambda)$ is defined by the generating relation \cite{Zhang}:

$$\sum\limits_{k=0}^{\infty} \mathcal{M}_{k}(n,\lambda)\frac{t^k}{k!}=\frac{1-\lambda (-e^t)^{(n+1)}}{\lambda e^{t}+1}.\eqno(1.20)$$

We note that
$$\mathcal{M}_{k}(n,1)=M_{k}(n).\eqno(1.21)$$

Also, for even $n$, we have
$$\mathcal{S}_k(n,-\lambda)=\mathcal{M}_k(n,\lambda).\eqno(1.22)$$

Motivated by the above mentioned work on symmetry identities for 1-variable Apostol polynomials and due to the importance of the 2-variable forms of the special polynomials in this paper, we derived symmetric identities for the 2-variable Apostol type polynomials. Symmetric identities for certain members belonging to this family are also obtained.\\
\vspace{.35cm}

\noindent
{\large{\bf 2. Symmetric identities for 2-variable Apostol type polynomials}}\\
\vspace{.35cm}

In this section, we derive several symmetric identities for the 2-variable Apostol type polynomials. For this we prove the following theorems:\\

\noindent
{\bf Theorem 2.1.}~{\em For all integers $c>0,~d >0,~n \geq0, ~m \in \mathbb{N}, ~\lambda,\mu,\nu \in \mathbb{C}$, the following symmetry identity for the 2-variable Apostol type polynomials $_p\mathcal{F}_n^{(m)}(x,y;\lambda;\mu,\nu)$ holds true}:
\[\begin{split}
&\sum_{k=0}^{n}{n\choose k}c^{n-k}d^{\nu+k} {}_p\mathcal{F}_{n-k}^{(m)}(dx,dy;\lambda;\mu,\nu)\sum_{l=0}^{k}{k\choose l}\mathcal{S}_{l}(c-1;-\lambda){}_p\mathcal{F}_{k-l}^{(m-1)}(cX,cY;\lambda;\mu,\nu)\\
&=\sum_{k=0}^{n}{n\choose k}d^{n-k}c^{\nu+k} {}_p\mathcal{F}_{n-k}^{(m)}(cx,cy;\lambda;\mu,\nu)\sum_{l=0}^{k}{k\choose l}\mathcal{S}_{l}(d-1;-\lambda){}_p\mathcal{F}_{k-l}^{(m-1)}
(dX,dY;\lambda;\mu,\nu).\hspace{.2in}(2.1)\\
\end{split}\]

\noindent
{\bf Proof.}~Let

$$G(t):= \frac{2^{\mu(2m-1)}t^{\nu(2m-1)}e^{cdxt}\phi(y,cdt)(\lambda e^{cdt}+1)e^{cdXt}\phi(Y,cdt)}{(\lambda e^{ct}+1)^m (\lambda e^{dt}+1)^m}.\eqno(2.2)$$

Since the expression for $G(t)$ is symmetric in $c$ and $d$, we can expand $G(t)$ into series in two ways, to obtain
\[\begin{split}
G(t)&=\frac{1}{c^{\nu m} d^{\nu (m-1)}}\left(\frac{2^\mu c^\nu t^\nu}{\lambda e^{ct}+1}\right)^m
e^{cdxt}\phi(y,cdt)\left(\frac{2^\mu d^\nu t^\nu}{\lambda e^{dt}+1}\right)^{m-1}\left(\frac{\lambda e^{cdt}+1}{\lambda e^{dt}+1}\right)e^{cdXt}\phi(Y,cdt)\\
&=\frac{1}{c^{\nu m} d^{\nu (m-1)}}\left(\frac{2^\mu (ct)^\nu}{\lambda e^{ct}+1}\right)^m
e^{cdxt}\phi(dy,ct)\left(\frac{\lambda e^{cdt}+1}{\lambda e^{dt}+1}\right)\left(\frac{2^\mu (dt)^\nu}{\lambda e^{dt}+1}\right)^{m-1}e^{cdXt}\phi(cY,dt)\\
&=\frac{1}{c^{\nu m} d^{\nu (m-1)}}\Big(\sum_{n=0}^{\infty}{}_p\mathcal{F}_{n}^{(m)}(dx,dy;\lambda;\mu,\nu)\frac{(ct)^n}{n!}\Big)\Big(\sum_{l=0}^{\infty}\mathcal{S}_{l}(c-1;-\lambda)\frac{(dt)^l}{l!}
\Big)\\
&\Big(\sum_{k=0}^{\infty}{}_p\mathcal{F}_{k}^{(m-1)}(cX,cY;\lambda;\mu,\nu)\frac{(dt)^k}{k!}\Big),\\
\end{split}\]
which on using \cite[p.890 Corollary 2]{Srivastava1} gives
\[\begin{split}
G(t)=&\frac{1}{c^{\nu m} d^{\nu m}}\sum_{n=0}^{\infty}\Big[\sum_{k=0}^{n}{n\choose k}c^{n-k}d^{\nu+k}{}_p\mathcal{F}_{n-k}^{(m)}(dx,dy;\lambda;\mu,\nu)\sum_{l=0}^{k}{k\choose l}\mathcal{S}_{l}(c-1;-\lambda)\\
&{}_p\mathcal{F}_{k-l}^{(m-1)}(cX,cY;\lambda;\mu,\nu)\Big]\frac{t^n}{n!}.\hspace{3.4in}(2.3)\\
\end{split}\]
In a similar manner, we have

\[\begin{split}
G(t)&=\frac{1}{d^{\nu m} c^{\nu (m-1)}}\left(\frac{2^\mu (dt)^\nu}{\lambda e^{dt}+1}\right)^m
e^{cdxt}\phi(cy,dt)\left(\frac{\lambda e^{cdt}+1}{\lambda e^{ct}+1}\right)\left(\frac{2^\mu (ct)^\nu}{\lambda e^{ct}+1}\right)^{m-1}e^{cdXt}\phi(dY,ct)\\
&=\frac{1}{d^{\nu m} c^{\nu (m-1)}}\Big(\sum_{n=0}^{\infty}{}_p\mathcal{F}_{n}^{(m)}(cx,cy;\lambda;\mu,\nu)\frac{(dt)^n}{n!}\Big)\Big(\sum_{l=0}^{\infty}\mathcal{S}_{l}(d-1;-\lambda)\frac{(ct)^l}{l!}
\Big)\\
&\Big(\sum_{k=0}^{\infty}{}_p\mathcal{F}_{k}^{(m-1)}(dX,dY;\lambda;\mu,\nu)\frac{(ct)^k}{k!}\Big)\\
&=\frac{1}{d^{\nu m} c^{\nu m}}\sum_{n=0}^{\infty}\Big[\sum_{k=0}^{n}{n\choose k}d^{n-k}c^{\nu+k}{}_p\mathcal{F}_{n-k}^{(m)}(cx,cy;\lambda;\mu,\nu)\sum_{l=0}^{k}{k\choose l}\mathcal{S}_{l}(d-1;-\lambda)\\
&{}_p\mathcal{F}_{k-l}^{(m-1)}(dX,dY;\lambda;\mu,\nu)\Big]\frac{t^n}{n!}.\hspace{3.8in}(2.4)\\
\end{split}\]

Equating the coefficients of same power of $t$ in r.h.s. of equations (2.3) and (2.4), we are led to assertion (2.1).\\

\noindent
{\bf Theorem 2.2.}~{\em For each pair of positive integers $c,d$ and for all integers $n \geq0, ~m \in \mathbb{N}, ~\lambda,\mu,\nu \in \mathbb{C}$, the following symmetry identity for the 2-variable Apostol type polynomials $_p\mathcal{F}_n^{(m)}(x,y;\lambda;\mu,\nu)$ holds true}:
\[\begin{split}
&\sum\limits_{k=0}^n {n \choose k} \sum\limits_{i=0}^{c-1}\sum\limits_{j=0}^{d-1} (-1)^{i+j}\lambda^{i+j} c^k d^{n-k} {}_p\mathcal{F}_n^{(m)}\Big(dx+\frac{d}{c}i,dy;\lambda;\mu,\nu\Big){}_p\mathcal{F}_n^{(m)}\Big(cX+\frac{c}{d}j,cY;\lambda;\mu,\nu\Big)\\
=&\sum\limits_{k=0}^n {n \choose k} \sum\limits_{i=0}^{d-1}\sum\limits_{j=0}^{c-1}(-1)^{i+j} \lambda^{i+j} d^k c^{n-k} {}_p\mathcal{F}_n^{(m)}\Big(cx+\frac{c}{d}i,cy;\lambda;\mu,\nu\Big){}_p\mathcal{F}_n^{(m)}\Big(dX+\frac{d}{c}j,dY;\lambda;\mu,\nu\Big).\hspace{.05in}(2.5)\\
\end{split}\]

\noindent
{\bf Proof.}~Let
$$H(t):= \frac{2^{2\mu m}t^{2\nu m}e^{cdxt}\phi(y,cdt)({\lambda}^c e^{cdt}+1)({\lambda}^d e^{cdt}+1) e^{cdXt}\phi(Y,cdt)}{(\lambda e^{ct}+1)^{m+1} (\lambda e^{dt}+1)^{m+1}}.\eqno(2.6)$$

Since the expression for $H(t)$ is symmetric in $c$ and $d$, we can expand $H(t)$ into series in two ways, to obtain
\[\begin{split}
H(t)=&\frac{1}{c^{\nu m}d^{\nu m}}\left(\frac{2^\mu c^\nu t^\nu}{\lambda e^{ct}+1}\right)^m
e^{cdxt}\phi(y,cdt)\Big(\frac{\lambda^c e^{cdt}+1}{\lambda e^{dt}+1}\Big)\left(\frac{2^\mu d^\nu t^\nu}{\lambda e^{dt}+1}\right)^{m}e^{cdXt}\phi(Y,cdt)\Big(\frac{\lambda^d e^{cdt}+1}{\lambda e^{ct}+1}\Big)\\
&=\frac{1}{c^{\nu m}d^{\nu m}}\left(\frac{2^\mu {(ct)}^\nu}{\lambda e^{ct}+1}\right)^m
e^{cdxt}\phi(dy,ct)\sum\limits_{i=0}^{c-1}(-\lambda)^i e^{dti}\left(\frac{2^\mu {(dt)}^\nu}{\lambda e^{dt}+1}\right)^{m}e^{cdXt}\phi(cY,dt)\sum\limits_{j=0}^{d-1}(-\lambda)^j e^{ctj}\\
&=\frac{1}{c^{\nu m}d^{\nu m}}\sum\limits_{i=0}^{c-1}(-\lambda)^i\left(\frac{2^\mu {(ct)}^\nu}{\lambda e^{ct}+1}\right)^m
e^{\big(dx+\frac{d}{c}i\big)ct}\phi(dy,ct)\sum\limits_{j=0}^{d-1}(-\lambda)^j \left(\frac{2^\mu {(dt)}^\nu}{\lambda e^{dt}+1}\right)^{m} e^{\big(cX+\frac{c}{d}j\big)dt}\phi(cY,dt) \\
&=\frac{1}{c^{\nu m}d^{\nu m}}\left(\sum\limits_{i=0}^{c-1}(-\lambda)^i \sum\limits_{n=0}^\infty {}_p\mathcal{F}_{n}^{(m)}\left(dx+\frac{d}{c}i,dy;\lambda;\mu,\nu\right)\frac{(ct)^n}{n!}\right)\\
\end{split}\]
\[\begin{split}
& \left(\sum\limits_{j=0}^{d-1}(-\lambda)^j \sum\limits_{n=0}^\infty {}_p\mathcal{F}_{n}^{(m)}\left(cX+\frac{c}{d}j,cY;\lambda;\mu,\nu\right)\frac{(dt)^n}{n!}\right)\\
&=\frac{1}{c^{\nu m}d^{\nu m}} \sum\limits_{n=0}^\infty \sum\limits_{k=0}^n {n \choose k} \sum\limits_{i=0}^{c-1}
\sum\limits_{j=0}^{d-1}(-\lambda)^{i+j}c^k d^{n-k}{}_p\mathcal{F}_{k}^{(m)}\left(dx+\frac{d}{c}i,dy;\lambda;\mu,\nu\right)\\
&{}_p\mathcal{F}_{n-k}^{(m)}\left(cX+\frac{c}{d}j,cY;\lambda;\mu,\nu\right)\frac{t^n}{n!}.\hspace{3.6in}(2.7)\\
\end{split}\]
In a similar manner, we have
\[\begin{split}
H(t)&=\frac{1}{c^{\nu m}d^{\nu m}} \sum\limits_{n=0}^\infty \sum\limits_{k=0}^n {n \choose k} \sum\limits_{i=0}^{d-1}
\sum\limits_{j=0}^{c-1}(-\lambda)^{i+j}d^k c^{n-k}{}_p\mathcal{F}_{k}^{(m)}\left(cx+\frac{c}{d}i,cy;\lambda;\mu,\nu\right)\\
&{}_p\mathcal{F}_{n-k}^{(m)}\left(dX+\frac{d}{c}j,dY;\lambda;\mu,\nu\right)\frac{t^n}{n!}.\hspace{3.2in}(2.8)\\
\end{split}\]

Equating the coefficients of same power of $t$ in r.h.s. of equations (2.7) and (2.8), we are led to assertion (2.2).\\

By taking suitable values of the parameters in equation (2.1) and (2.5) and in view of relations (1.12)-(1.14), we can find the symmetric identities for other mixed special polynomials related to $_p\mathcal{F}_n^{(m)}(x,y;\lambda;\mu,\nu)$. We present the symmetric identities for these polynomials in the following table:\\

\noindent
{\bf Table~2.1.~Symmetry identities for the special cases of the 2VATP $_p\mathcal{F}_n^{(m)}(x,y;\lambda;\mu,\nu)$.}\\
\\
{\tiny{
\begin{tabular}{|l|l|l|l|l|}
\hline
&&&&\\
{\bf S.} & {\bf Values of}  & {\bf Relation between} & {\bf Name of the} &  {\bf Symmetry identities  of the resultant special polynomials}\\

 {\bf No.} &  {\bf the para-}  &  {\bf the 2VATP}  & {\bf resultant} & \\
   & {\bf meters}  &  $_p\mathcal{F}_n^{(m)}(x,y;\lambda;\mu,\nu)$   &{\bf special}   &  \\
  &    &  {\bf and its special case}  &   {\bf polynomials}  &   \\
\hline

{\bf I.}  &  $\lambda \rightarrow -\lambda$  & $(-1)^m $  & 2-variable& $\sum\limits_{k=0}^{n}{n\choose k}c^{n-k}d^{\nu+k} {}_p\mathfrak{B}_{n-k}^{(m)}(dx,dy;\lambda)\sum\limits_{l=0}^{k}{k\choose l}\mathcal{S}_{l}(c-1;\lambda){}_p\mathfrak{B}_{k-l}^{(m-1)}(cX,cY;\lambda)$\\

  &  $\mu=0$   & ${}_p\mathcal{F}_n^{(m)}(x,y;-\lambda;0,1)$   &   Apostol & $=\sum\limits_{k=0}^{n}{n\choose k}d^{n-k}c^{\nu+k} {}_p\mathfrak{B}_{n-k}^{(m)}(cx,cy;\lambda)\sum\limits_{l=0}^{k}{k\choose l}\mathcal{S}_{l}(d-1;\lambda){}_p\mathfrak{B}_{k-l}^{(m-1)}(dX,dY;\lambda)$  \\

     &   $\nu=1$  &  $={}_p\mathfrak{B}_n^{(m)}(x,y;\lambda)$ &   Bernoulli &\\

  &     &     &   polynomials    & $\sum\limits_{k=0}^n {n \choose k} \sum\limits_{i=0}^{c-1}\sum\limits_{j=0}^{d-1} (-1)^{i+j}(-\lambda)^{i+j} c^k d^{n-k} {}_p\mathfrak{B}_n^{(m)}\Big(dx+\frac{d}{c}i,dy;\lambda){}_p\mathfrak{B}_n^{(m)}\Big(cX+\frac{c}{d}j,cY;\lambda\Big)$ \\

     &    &   & (2VABP)  &  $\sum\limits_{k=0}^n {n \choose k} \sum\limits_{i=0}^{d-1}\sum\limits_{j=0}^{c-1}(-1)^{i+j} (-\lambda)^{i+j} d^k c^{n-k} {}_p\mathfrak{B}_n^{(m)}\Big(cx+\frac{c}{d}i,cy;\lambda\Big){}_p\mathfrak{B}_n^{(m)}\Big(dX+\frac{d}{c}j,dY;\lambda\Big)$  \\
\hline

{\bf II.}  &  $\mu=1$  & $ {}_p\mathcal{F}_n^{(m)}(x,y;\lambda;1,0)$  & 2-variable& $\sum\limits_{k=0}^{n}{n\choose k}c^{n-k}d^{\nu+k} {}_p\mathfrak{E}_{n-k}^{(m)}(dx,dy;\lambda)\sum\limits_{l=0}^{k}{k\choose l}\mathcal{M}_{l}(c-1;\lambda){}_p\mathfrak{E}_{k-l}^{(m-1)}(cX,cY;\lambda)$\\

  &   $\nu=0$   &  $={}_p\mathfrak{E}_n^{(m)}(x,y;\lambda)$  &   Apostol & $=\sum\limits_{k=0}^{n}{n\choose k}d^{n-k}c^{\nu+k} {}_p\mathfrak{E}_{n-k}^{(m)}(cx,cy;\lambda)\sum\limits_{l=0}^{k}{k\choose l}\mathcal{M}_{l}(d-1;\lambda){}_p\mathfrak{E}_{k-l}^{(m-1)}(dX,dY;\lambda)$  \\

     &    &   &   Euler &\\

  &     &     &   polynomials    & $\sum\limits_{k=0}^n {n \choose k} \sum\limits_{i=0}^{c-1}\sum\limits_{j=0}^{d-1} (-1)^{i+j}\lambda^{i+j} c^k d^{n-k} {}_p\mathfrak{E}_n^{(m)}\Big(dx+\frac{d}{c}i,dy;\lambda){}_p\mathfrak{E}_n^{(m)}\Big(cX+\frac{c}{d}j,cY;\lambda\Big)$ \\

     &    &   & (2VAEP)  &  $\sum\limits_{k=0}^n {n \choose k} \sum\limits_{i=0}^{d-1}\sum\limits_{j=0}^{c-1}(-1)^{i+j} \lambda^{i+j} d^k c^{n-k} {}_p\mathfrak{E}_n^{(m)}\Big(cx+\frac{c}{d}i,cy;\lambda\Big){}_p\mathfrak{E}_n^{(m)}\Big(dX+\frac{d}{c}j,dY;\lambda\Big)$  \\
\hline
{\bf III.}  &  $\mu=1$  & $ {}_p\mathcal{F}_n^{(m)}(x,y;\lambda;1,1)$  & 2-variable& $\sum\limits_{k=0}^{n}{n\choose k}c^{n-k}d^{\nu+k} {}_p\mathcal{G}_{n-k}^{(m)}(dx,dy;\lambda)\sum\limits_{l=0}^{k}{k\choose l}\mathcal{M}_{l}(c-1;\lambda){}_p\mathcal{G}_{k-l}^{(m-1)}(cX,cY;\lambda)$\\

  &   $\nu=1$   &  $={}_p\mathcal{G}_n^{(m)}(x,y;\lambda)$  &   Apostol & $=\sum\limits_{k=0}^{n}{n\choose k}d^{n-k}c^{\nu+k} {}_p\mathcal{G}_{n-k}^{(m)}(cx,cy;\lambda)\sum\limits_{l=0}^{k}{k\choose l}\mathcal{M}_{l}(d-1;\lambda){}_p\mathcal{G}_{k-l}^{(m-1)}(dX,dY;\lambda)$  \\

     &    &   &   Genocchi &\\

  &     &     &   polynomials    & $\sum\limits_{k=0}^n {n \choose k} \sum\limits_{i=0}^{c-1}\sum\limits_{j=0}^{d-1} (-1)^{i+j}\lambda^{i+j} c^k d^{n-k} {}_p\mathcal{G}_n^{(m)}\Big(dx+\frac{d}{c}i,dy;\lambda){}_p\mathcal{G}_n^{(m)}\Big(cX+\frac{c}{d}j,cY;\lambda\Big)$ \\

     &    &   & (2VAGP)  &  $\sum\limits_{k=0}^n {n \choose k} \sum\limits_{i=0}^{d-1}\sum\limits_{j=0}^{c-1}(-1)^{i+j} \lambda^{i+j} d^k c^{n-k} {}_p\mathcal{G}_n^{(m)}\Big(cx+\frac{c}{d}i,cy;\lambda\Big){}_p\mathcal{G}_n^{(m)}\Big(dX+\frac{d}{c}j,dY;\lambda\Big)$  \\
\hline

\end{tabular}}}\\
\vspace{.35cm}

\noindent
{\bf Note.}~We note that for $\lambda=1$, the results derived in Table 2.1 for the 2VABP $_p\mathfrak{B}_n^{(m)}(x,y;\lambda)$, 2VAEP $_p\mathfrak{E}_n^{(m)}(x,y;\lambda)$ and 2VAGP $_p\mathcal{G}_n^{(m)}(x,y;\lambda)$ give the corresponding results for the 2-variable Bernoulli polynomials (2VBP) (of order $m$)  $_pB_n^{(m)}(x,y)$, 2-variable Euler polynomials (2VEP) (of order $m$) $_pE_n^{(m)}(x,y)$ and 2-variable Genocchi polynomials (2VGP) (of order $m$) $_pG_n^{(m)}(x,y)$, respectively. We present these results in the Table 2.2:\\
\vspace{.35cm}

\noindent
{\bf Table~2.2.~Symmetry identities for the $_pB_n^{(m)}(x,y)$, $_pE_n^{(m)}(x,y)$ and $_pG_n^{(m)}(x,y)$.}\\
\\
{\tiny{
\begin{tabular}{|l|l|l|l|l|}
\hline
&&&\\
{\bf S.}  & {\bf Relation between} & {\bf Name of the} &  {\bf Symmetry identities  of the resultant special polynomials}\\

 {\bf No.}   &  {\bf the special}  & {\bf resultant} & \\
   & {\bf  cases}   &{\bf special}   &  \\
     &    &   {\bf polynomials}  &   \\
\hline

{\bf I.}   & ${}_p\mathfrak{B}_n^{(m)}(x,y;1)$  & 2-variable& $\sum\limits_{k=0}^{n}{n\choose k}c^{n-k}d^{\nu+k} {}_pB_{n-k}^{(m)}(dx,dy)\sum\limits_{l=0}^{k}{k\choose l}S_{l}(c-1){}_pB_{k-l}^{(m-1)}(cX,cY)$\\
  &  $={}_pB_n^{(m)}(x,y)$  &  Bernoulli   & $=\sum\limits_{k=0}^{n}{n\choose k}d^{n-k}c^{\nu+k} {}_pB_{n-k}^{(m)}(cx,cy)\sum\limits_{l=0}^{k}{k\choose l}S_{l}(d-1){}_pB_{k-l}^{(m-1)}(dX,dY)$  \\

   &   & polynomials & \\

      &     &  (2VBP)     & $\sum\limits_{k=0}^n {n \choose k} \sum\limits_{i=0}^{c-1}\sum\limits_{j=0}^{d-1}  c^k d^{n-k} {}_pB_n^{(m)}\Big(dx+\frac{d}{c}i,dy){}_pB_n^{(m)}\Big(cX+\frac{c}{d}j,cY\Big)$ \\

         &   &  of order $m$ &  $=\sum\limits_{k=0}^n {n \choose k} \sum\limits_{i=0}^{d-1}\sum\limits_{j=0}^{c-1} d^k c^{n-k} {}_pB_n^{(m)}\Big(cx+\frac{c}{d}i,cy\Big){}_pB_n^{(m)}\Big(dX+\frac{d}{c}j,dY\Big)$  \\
\hline

{\bf II.}   & ${}_p\mathfrak{E}_n^{(m)}(x,y;1)$  & 2-variable& $\sum\limits_{k=0}^{n}{n\choose k}c^{n-k}d^{\nu+k} {}_pE_{n-k}^{(m)}(dx,dy)\sum\limits_{l=0}^{k}{k\choose l}M_{l}(c-1){}_pE_{k-l}^{(m-1)}(cX,cY)$\\
  &  $={}_pE_n^{(m)}(x,y)$  &  Euler   & $=\sum\limits_{k=0}^{n}{n\choose k}d^{n-k}c^{\nu+k} {}_pE_{n-k}^{(m)}(cx,cy)\sum\limits_{l=0}^{k}{k\choose l}M_{l}(d-1){}_pE_{k-l}^{(m-1)}(dX,dY)$  \\

   &   & polynomials & \\

      &     &  (2VEP)     & $\sum\limits_{k=0}^n {n \choose k} \sum\limits_{i=0}^{c-1}\sum\limits_{j=0}^{d-1}  c^k d^{n-k} {}_pE_n^{(m)}\Big(dx+\frac{d}{c}i,dy){}_pE_n^{(m)}\Big(cX+\frac{c}{d}j,cY\Big)$ \\

         &   &  of order $m$  &  $=\sum\limits_{k=0}^n {n \choose k} \sum\limits_{i=0}^{d-1}\sum\limits_{j=0}^{c-1} d^k c^{n-k} {}_pE_n^{(m)}\Big(cx+\frac{c}{d}i,cy\Big){}_pE_n^{(m)}\Big(dX+\frac{d}{c}j,dY\Big)$  \\
\hline

{\bf III.}   & ${}_p\mathcal{G}_n^{(m)}(x,y;1)$  & 2-variable& $\sum\limits_{k=0}^{n}{n\choose k}c^{n-k}d^{\nu+k} {}_pG_{n-k}^{(m)}(dx,dy)\sum\limits_{l=0}^{k}{k\choose l}M_{l}(c-1){}_pG_{k-l}^{(m-1)}(cX,cY)$\\
  &  $={}_pG_n^{(m)}(x,y)$  &  Genocchi   & $=\sum\limits_{k=0}^{n}{n\choose k}d^{n-k}c^{\nu+k} {}_pG_{n-k}^{(m)}(cx,cy)\sum\limits_{l=0}^{k}{k\choose l}M_{l}(d-1){}_pG_{k-l}^{(m-1)}(dX,dY)$  \\

   &   & polynomials & \\

      &     &  (2VGP)     & $\sum\limits_{k=0}^n {n \choose k} \sum\limits_{i=0}^{c-1}\sum\limits_{j=0}^{d-1}  c^k d^{n-k} {}_pG_n^{(m)}\Big(dx+\frac{d}{c}i,dy){}_pG_n^{(m)}\Big(cX+\frac{c}{d}j,cY\Big)$ \\

         &   & of order $m$  &  $=\sum\limits_{k=0}^n {n \choose k} \sum\limits_{i=0}^{d-1}\sum\limits_{j=0}^{c-1} d^k c^{n-k} {}_pG_n^{(m)}\Big(cx+\frac{c}{d}i,cy\Big){}_pG_n^{(m)}\Big(dX+\frac{d}{c}j,dY\Big)$  \\
\hline
\end{tabular}}}\\
\vspace{.35cm}

\noindent
{\bf Note.}~We note that for $m=1$, the results derived above for the 2-variable Bernoulli polynomials (2VBP) (of order $m$) $_pB_n^{(m)}(x,y)$, 2-variable Euler polynomials (2VEP) (of order $m$) $_pE_n^{(m)}(x,y)$ and 2-variable Genocchi polynomials (2VGP)(of order $m$) $_pG_n^{(m)}(x,y)$ give the corresponding results for the 2-variable Bernoulli polynomials (2VBP) $_pB_n(x,y)$, 2-variable Euler polynomials (2VEP) $_pE_n(x,y)$ and 2-variable Genocchi polynomials (2VGP) $_pG_n(x,y)$.\\

In the next section, symmetry identities for certain members belonging to the 2-variable Apostol type polynomials are obtained.\\
\vspace{.35cm}

\noindent
{\large{\bf 3.~Symmetry identities for certain members belonging to 2-variable\\ Apostol type polynomials}}
\vspace{.35cm}

Certain members belonging to the 2VATP $_p\mathcal{F}_n^{(m)}(x,y;\lambda;\mu,\nu)$ family are considered in \cite{ATP}. These special members are obtained by making suitable choice for the function $\phi(y,t)$ in equation (1.15).\\

We recall that by taking $\phi(y,t)=e^{yt^s}$ (for which the 2VGP $p_n(x,y)$ reduce to the Gould-Hopper polynomials (GHP) $H_n^{(s)}(x,y))$ in the l.h.s. of generating function (1.15), we get the following generating function for the Gould-Hopper Apostol type polynomials ${}_{H^{(s)}}\mathcal{F}_n^{(m)}(x,y;\lambda;\mu,\nu)$ \cite{ATP}:

$$\Big(\frac{2^\mu~ t^\nu}{\lambda e^{t}+1}\Big)^m e^{xt+yt^s}=\sum\limits_{n=0}^{\infty} {}_{H^{(s)}}\mathcal{F}_n^{(m)}(x,y;\lambda;\mu,\nu)\frac{t^{n}}{n!}.\eqno(3.1)$$

Consequently, by taking $\phi(y,cdt)=e^{y(cdt)^s}$ in equation (2.2) and proceeding on the same lines as in Theorem 2.1 and Theorem 2.2, we get the following consequences of Theorem 2.1 and Theorem 2.2.\\

\noindent
{\bf Corollary 3.1.}~{\em For all integers $c>0,~d >0,~n \geq0, ~m \in \mathbb{N}, ~\lambda,\mu,\nu \in \mathbb{C}$, the following symmetry identity for the Gould-Hopper Apostol type polynomials $_{H^{(s)}}\mathcal{F}_n^{(m)}(x,y;\lambda;\mu,\nu)$ holds true:}
\[\begin{split}
&\sum_{k=0}^{n}{n\choose k}c^{n-k}d^{\nu+k} {}_{H^{(s)}}\mathcal{F}_{n-k}^{(m)}(dx,d^{s}y;\lambda;\mu,\nu)\sum_{l=0}^{k}{k\choose l}\mathcal{S}_{l}(c-1;-\lambda){}_{H^{(s)}}\mathcal{F}_{k-l}^{(m-1)}(cX,c^{s}Y;\lambda;\mu,\nu)\\
&=\sum_{k=0}^{n}{n\choose k}d^{n-k}c^{\nu+k} {}_{H^{(s)}}\mathcal{F}_{n-k}^{(m)}(cx,c^{s}y;\lambda;\mu,\nu)\sum_{l=0}^{k}{k\choose l}\mathcal{S}_{l}(d-1;-\lambda){}_{H^{(s)}}\mathcal{F}_{k-l}^{(m-1)}
(dX,d^{s}Y;\lambda;\mu,\nu).\\
&\hspace{5.24in}(3.2)
\end{split}\]

\noindent
{\bf Corollary 3.2.}~{\em For each pair of positive integers $c, ~d$ and for all integers $n \geq0, ~m \in \mathbb{N}, ~\lambda,\mu,\nu \in \mathbb{C}$, the following symmetry identity for the Gould-Hopper Apostol type polynomials $_{H^{(s)}}\mathcal{F}_n^{(m)}(x,y;\lambda;\mu,\nu)$ holds true}:
\[\begin{split}
&\sum\limits_{k=0}^n {n \choose k} \sum\limits_{i=0}^{c-1}\sum\limits_{j=0}^{d-1} (-1)^{i+j}\lambda^{i+j} c^k d^{n-k} {}_{H^{(s)}}\mathcal{F}_n^{(m)}\Big(dx+\frac{d}{c}i,d^{s}y;\lambda;\mu,\nu\Big){}_{H^{(s)}}\mathcal{F}_n^{(m)}\Big(cX+\frac{c}{d}j,c^{s}Y;\lambda;\mu,\nu\Big)\\
&=\sum\limits_{k=0}^n {n \choose k} \sum\limits_{i=0}^{d-1}\sum\limits_{j=0}^{c-1}(-1)^{i+j} \lambda^{i+j} d^k c^{n-k} {}_{H^{(s)}}\mathcal{F}_n^{(m)}\Big(cx+\frac{c}{d}i,c^{s}y;\lambda;\mu,\nu\Big){}_{H^{(s)}}\mathcal{F}_n^{(m)}\Big(dX+\frac{d}{c}j,d^{s}Y;\lambda;\mu,\nu\Big).\\
&\hspace{5.20in}(3.3)
\end{split}\]

By taking suitable values of the parameters in equation (3.2) and (3.3) and in view of relations (1.12)-(1.14), we can find the symmetric identities for other mixed special polynomials related to $_{H^{(s)}}\mathcal{F}_n^{(m)}(x,y;\lambda;\mu,\nu)$. We present the symmetric identities for these polynomials in the following table:\\
\vspace{.35cm}

\noindent
{\bf Table~3.1.~Symmetry identities for the Special cases of the GHATP $_{H^{(s)}}\mathcal{F}_n^{(m)}(x,y;\lambda;\mu,\nu)$}\\
\\
{\tiny{
\begin{tabular}{|l|l|l|l|l|}
\hline
&&&\\
{\bf S.} & {\bf Values of}  & {\bf Relation between}  &  {\bf Name and Symmetry identities  of the resultant special polynomials}\\

 {\bf No.} &  {\bf the para-}  &  {\bf the GHATP}   & \\
   & {\bf meters}  &  $_{H^{(s)}}\mathcal{F}_n^{(m)}(x,y;\lambda;\mu,\nu)$      &  \\
  &    &  {\bf and its special case}    &   \\
\hline
&&&\\
{\bf I.}  &  $\lambda \rightarrow -\lambda$  & $(-1)^m $  & Gould-Hopper Apostol Bernoulli polynomials \\

  &  $\mu=0$   & ${}_{H^{(s)}}\mathcal{F}_n^{(m)}(x,y;-\lambda;0,1)$   &  $\sum\limits_{k=0}^{n}{n\choose k}c^{n-k}d^{\nu+k} {}_{H^{(s)}}\mathfrak{B}_{n-k}^{(m)}(dx,d^{s}y;\lambda)\sum\limits_{l=0}^{k}{k\choose l}\mathcal{S}_{l}(c-1;\lambda){}_{H^{(s)}}\mathfrak{B}_{k-l}^{(m-1)}(cX,c^{s}Y;\lambda)$ \\

     &   $\nu=1$  &  $={}_{H^{(s)}}\mathfrak{B}_n^{(m)}(x,y;\lambda)$  & $=\sum\limits_{k=0}^{n}{n\choose k}d^{n-k}c^{\nu+k} {}_{H^{(s)}}\mathfrak{B}_{n-k}^{(m)}(cx,c^{s}y;\lambda)\sum\limits_{l=0}^{k}{k\choose l}\mathcal{S}_{l}(d-1;\lambda){}_{H^{(s)}}\mathfrak{B}_{k-l}^{(m-1)}(dX,d^{s}Y;\lambda)$\\
&&&\\
  &     &         & $\sum\limits_{k=0}^n {n \choose k} \sum\limits_{i=0}^{c-1}\sum\limits_{j=0}^{d-1} (-1)^{i+j}(-\lambda)^{i+j} c^k d^{n-k} {}_{H^{(s)}}\mathfrak{B}_n^{(m)}\Big(dx+\frac{d}{c}i,d^{s}y;\lambda){}_{H^{(s)}}\mathfrak{B}_n^{(m)}\Big(cX+\frac{c}{d}j,c^{s}Y;\lambda\Big)$ \\

     &    &     &  $\sum\limits_{k=0}^n {n \choose k} \sum\limits_{i=0}^{d-1}\sum\limits_{j=0}^{c-1}(-1)^{i+j} (-\lambda)^{i+j} d^k c^{n-k} {}_{H^{(s)}}\mathfrak{B}_n^{(m)}\Big(cx+\frac{c}{d}i,c^{s}y;\lambda\Big){}_{H^{(s)}}\mathfrak{B}_n^{(m)}\Big(dX+\frac{d}{c}j,d^{s}Y;\lambda\Big)$  \\
\hline

{\bf II.}  &  $\mu=1$  & $ {}_{H^{(s)}}\mathcal{F}_n^{(m)}(x,y;\lambda;1,0)$  & Gould-Hopper Apostol Euler polynomials\\

  &   $\nu=0$   &  $={}_{H^{(s)}}\mathfrak{E}_n^{(m)}(x,y;\lambda)$   & $\sum\limits_{k=0}^{n}{n\choose k}c^{n-k}d^{\nu+k} {}_{H^{(s)}}\mathfrak{E}_{n-k}^{(m)}(dx,d^{s}y;\lambda)\sum\limits_{l=0}^{k}{k\choose l}\mathcal{M}_{l}(c-1;\lambda){}_{H^{(s)}}\mathfrak{E}_{k-l}^{(m-1)}(cX,c^{s}Y;\lambda)$\\

     &    &    &$=\sum\limits_{k=0}^{n}{n\choose k}d^{n-k}c^{\nu+k} {}_{H^{(s)}}\mathfrak{E}_{n-k}^{(m)}(cx,c^{s}y;\lambda)\sum\limits_{l=0}^{k}{k\choose l}\mathcal{M}_{l}(d-1;\lambda){}_{H^{(m)}}\mathfrak{E}_{k-l}^{(m-1)}(dX,d^{s}Y;\lambda)$  \\
&&&\\
  &     &        & $\sum\limits_{k=0}^n {n \choose k} \sum\limits_{i=0}^{c-1}\sum\limits_{j=0}^{d-1} (-1)^{i+j}\lambda^{i+j} c^k d^{n-k} {}_{H^{(s)}}\mathfrak{E}_n^{(m)}\Big(dx+\frac{d}{c}i,d^{s}y;\lambda){}_{H^{(s)}}\mathfrak{E}_n^{(m)}\Big(cX+\frac{c}{d}j,c^{s}Y;\lambda\Big)$ \\

     &    &     &  $\sum\limits_{k=0}^n {n \choose k} \sum\limits_{i=0}^{d-1}\sum\limits_{j=0}^{c-1}(-1)^{i+j} \lambda^{i+j} d^k c^{n-k} {}_{H^{(s)}}\mathfrak{E}_n^{(m)}\Big(cx+\frac{c}{d}i,c^{s}y;\lambda\Big){}_{H^{(s)}}\mathfrak{E}_n^{(m)}\Big(dX+\frac{d}{c}j,d^{s}Y;\lambda\Big)$  \\
\hline

{\bf III.}  &  $\mu=1$  & $ {}_{H^{(s)}}\mathcal{F}_n^{(m)}(x,y;\lambda;1,1)$  & Gould-Hopper Apostol Genocchi polynomials\\

  &   $\nu=1$   &  $={}_{H^{(s)}}\mathcal{G}_n^{(m)}(x,y;\lambda)$   & $\sum\limits_{k=0}^{n}{n\choose k}c^{n-k}d^{\nu+k} {}_{H^{(s)}}\mathcal{G}_{n-k}^{(m)}(dx,d^{s}y;\lambda)\sum\limits_{l=0}^{k}{k\choose l}\mathcal{M}_{l}(c-1;\lambda){}_{H^{(s)}}\mathcal{G}_{k-l}^{(m-1)}(cX,c^{s}Y;\lambda)$ \\

     &    &    &$=\sum\limits_{k=0}^{n}{n\choose k}d^{n-k}c^{\nu+k} {}_{H^{(s)}}\mathcal{G}_{n-k}^{(m)}(cx,c^{s}y;\lambda)\sum\limits_{l=0}^{k}{k\choose l}\mathcal{M}_{l}(d-1;\lambda){}_{H^{(s)}}\mathcal{G}_{k-l}^{(m-1)}(dX,d^{s}Y;\lambda)$ \\
&&&\\
  &     &        & $\sum\limits_{k=0}^n {n \choose k} \sum\limits_{i=0}^{c-1}\sum\limits_{j=0}^{d-1} (-1)^{i+j}\lambda^{i+j} c^k d^{n-k} {}_{H^{(s)}}\mathcal{G}_n^{(m)}\Big(dx+\frac{d}{c}i,d^{s}y;\lambda){}_{H^{(s)}}\mathcal{G}_n^{(m)}\Big(cX+\frac{c}{d}j,c^{s}Y;\lambda\Big)$ \\

     &    &     &  $\sum\limits_{k=0}^n {n \choose k} \sum\limits_{i=0}^{d-1}\sum\limits_{j=0}^{c-1}(-1)^{i+j} \lambda^{i+j} d^k c^{n-k} {}_{H^{(s)}}\mathcal{G}_n^{(m)}\Big(cx+\frac{c}{d}i,c^{s}y;\lambda\Big){}_{H^{(s)}}\mathcal{G}_n^{(m)}\Big(dX+\frac{d}{c}j,d^{s}Y;\lambda\Big)$  \\
\hline

\end{tabular}}}\\
\vspace{.35cm}

\noindent
{\bf Note.}~We note that for $\lambda=1$, the symmetry identities derived above for the GHABP $_{H^{(s)}}\mathfrak{B}_n^{(m)}(x,y;\lambda)$, GHAEP $_{H^{(s)}}\mathfrak{E}_n^{(m)}(x,y;\lambda)$ and GHAGP $_{H^{(s)}}\mathcal{G}_n^{(m)}(x,y;\lambda)$ give the corresponding identities for the Gould-Hopper Bernoulli polynomials (GHBP) (of order $m$)  $_{H^{(s)}}B_n^{(m)}(x,y)$, Gould-Hopper Euler polynomials (GHEP) (of order $m$) $_{H^{(s)}}E_n^{(m)}(x,y)$ and Gould-Hopper Genocchi polynomials (GHGP) (of order $m$) $_{H^{(s)}}G_n^{(m)}(x,y)$, respectively. Again, for $m=1$, we get the corresponding identities for the Gould-Hopper Bernoulli polynomials (GHBP) $_{H^{(s)}}B_n(x,y)$, Gould-Hopper Euler polynomials (GHEP) $_{H^{(s)}}E_n(x,y)$ and Gould-Hopper Genocchi polynomials (GHGP) $_{H^{(s)}}G_n(x,y)$.\\

Next, by taking $\phi(y,t)=C_0(-yt^s)$ (for which the 2VGP $p_n(x,y)$ reduce to the 2-variable generalized Laguerre polynomials (2VGLP) $_sL_n(y,x))$ in the l.h.s. of generating function (1.15), we get the following generating function for the 2-variable generalized Laguerre Apostol type polynomials ${}_{_sL}\mathcal{F}_n^{(m)}(x,y;\lambda;\mu,\nu)$ \cite{ATP}:

$$\Big(\frac{2^\mu~ t^\nu}{\lambda e^{t}+1}\Big)^m e^{xt}C_0(-yt^s)=\sum\limits_{n=0}^{\infty} {}_{_sL}\mathcal{F}_n^{(m)}(y,x;\lambda;\mu,\nu)\frac{t^{n}}{n!}.\eqno(3.4)$$

Consequently, by taking $\phi(y,cdt)=C_0(-y(cdt)^s)$ in equation (2.2) and proceeding on the same lines as in Theorem 2.1 and Theorem 2.2, we get the following consequences of Theorem 2.1 and Theorem 2.2.\\

\noindent
{\bf Corollary 3.3.}~{\em For all integers $c>0,~d >0,~n \geq0, ~m \in \mathbb{N}, ~\lambda,\mu,\nu \in \mathbb{C}$, the following symmetry identity for the 2-variable generalized Laguerre Apostol type polynomials $_{_sL}\mathcal{F}_n^{(m)}(x,y;\lambda;\mu,\nu)$ holds true:}
\[\begin{split}
&\sum_{k=0}^{n}{n\choose k}c^{n-k}d^{\nu+k} {}_{_sL}\mathcal{F}_{n-k}^{(m)}(d^{s}y,dx;\lambda;\mu,\nu)\sum_{l=0}^{k}{k\choose l}\mathcal{S}_{l}(c-1;-\lambda){}_{_sL}\mathcal{F}_{k-l}^{(m-1)}(c^{s}Y,cX;\lambda;\mu,\nu)\\
&=\sum_{k=0}^{n}{n\choose k}d^{n-k}c^{\nu+k} {}_{_sL}\mathcal{F}_{n-k}^{(m)}(c^{s}y,cx;\lambda;\mu,\nu)\sum_{l=0}^{k}{k\choose l}\mathcal{S}_{l}(d-1;-\lambda){}_{_sL}\mathcal{F}_{k-l}^{(m-1)}
(d^{s}Y,dX;\lambda;\mu,\nu).\\
&\hspace{5.14in}(3.5)
\end{split}\]

\noindent
{\bf Corollary 3.4.}~{\em For each pair of positive integers $c,~d$ and for all integers $n \geq0, ~m \in \mathbb{N}, ~\lambda,\mu,\nu \in \mathbb{C}$, the following symmetry identity for the 2-variable generalized Laguerre Apostol type polynomials $_{_sL}\mathcal{F}_n^{(m)}(x,y;\lambda;\mu,\nu)$ holds true}:
\[\begin{split}
&\sum\limits_{k=0}^n {n \choose k} \sum\limits_{i=0}^{c-1}\sum\limits_{j=0}^{d-1} (-1)^{i+j}\lambda^{i+j} c^k d^{n-k} {}_{_sL}\mathcal{F}_n^{(m)}\Big(d^{s}y,dx+\frac{d}{c}i;\lambda;\mu,\nu\Big){}_{_sL}\mathcal{F}_n^{(m)}\Big(c^{s}Y,cX+\frac{c}{d}j;\lambda;\mu,\nu\Big)\\
&=\sum\limits_{k=0}^n {n \choose k} \sum\limits_{i=0}^{d-1}\sum\limits_{j=0}^{c-1}(-1)^{i+j} \lambda^{i+j} d^k c^{n-k} {}_{_sL}\mathcal{F}_n^{(m)}\Big(c^{s}y,cx+\frac{c}{d}i;\lambda;\mu,\nu\Big){}_{_sL}\mathcal{F}_n^{(m)}\Big(d^{s}Y,dX+\frac{d}{c}j;\lambda;\mu,\nu\Big).\\
&\hspace{5.14in}(3.6)
\end{split}\]

By taking suitable values of the parameters in equation (3.2) and (3.3) and in view of relations (1.12)-(1.14), we can find the symmetric identities for other mixed special polynomials related to $_{_sL}\mathcal{F}_n^{(m)}(x,y;\lambda;\mu,\nu)$. We present the symmetric identities for these polynomials in the following table:\\
\vspace{.35cm}

\noindent
{\bf Table~3.2.~Symmetry identities for the Special cases of the 2VGLATP $_{_sL}\mathcal{F}_n^{(m)}(x,y;\lambda;\mu,\nu)$}\\
\\
{\tiny{
\begin{tabular}{|l|l|l|l|l|}
\hline
&&&\\
{\bf S.} & {\bf Values of}  & {\bf Relation between}  &  {\bf Name and Symmetry identities  of the resultant special polynomials}\\

 {\bf No.} &  {\bf the para-}  &  {\bf the 2VGLATP}   & \\
   & {\bf meters}  &  $_{_sL}\mathcal{F}_n^{(m)}(y,x;\lambda;\mu,\nu)$      &  \\
  &    &  {\bf and its special case}    &   \\
\hline
&&&\\
{\bf I.}  &  $\lambda \rightarrow -\lambda$  & $(-1)^m $  & 2-variable generalized Laguerre Apostol Bernoulli polynomials \\

  &  $\mu=0$   & ${}_{_sL}\mathcal{F}_n^{(m)}(y,x;-\lambda;0,1)$   &  $\sum\limits_{k=0}^{n}{n\choose k}c^{n-k}d^{\nu+k} {}_{_sL}\mathfrak{B}_{n-k}^{(m)}(d^{s}y,dx;\lambda)\sum\limits_{l=0}^{k}{k\choose l}\mathcal{S}_{l}(c-1;\lambda){}_{_sL}\mathfrak{B}_{k-l}^{(m-1)}(c^{s}Y,cX;\lambda)$ \\

     &   $\nu=1$  &  $={}_{_sL}\mathfrak{B}_n^{(m)}(y,x;\lambda)$  & $=\sum\limits_{k=0}^{n}{n\choose k}d^{n-k}c^{\nu+k} {}_{_sL}\mathfrak{B}_{n-k}^{(m)}(c^{s}y,cx;\lambda)\sum\limits_{l=0}^{k}{k\choose l}\mathcal{S}_{l}(d-1;\lambda){}_{_sL}\mathfrak{B}_{k-l}^{(m-1)}(d^{s}Y,dX;\lambda)$\\
&&&\\
  &     &         & $\sum\limits_{k=0}^n {n \choose k} \sum\limits_{i=0}^{c-1}\sum\limits_{j=0}^{d-1} (-1)^{i+j}(-\lambda)^{i+j} c^k d^{n-k} {}_{_sL}\mathfrak{B}_n^{(m)}\Big(d^{s}y,dx+\frac{d}{c}i;\lambda){}_{_sL}\mathfrak{B}_n^{(m)}\Big(c^{s}Y,cX+\frac{c}{d}j;\lambda\Big)$ \\

     &    &     &  $\sum\limits_{k=0}^n {n \choose k} \sum\limits_{i=0}^{d-1}\sum\limits_{j=0}^{c-1}(-1)^{i+j} (-\lambda)^{i+j} d^k c^{n-k} {}_{_sL}\mathfrak{B}_n^{(m)}\Big(c^{s}y,cx+\frac{c}{d}i;\lambda\Big){}_{_sL}\mathfrak{B}_n^{(m)}\Big(d^{s}Y,dX+\frac{d}{c}j;\lambda\Big)$  \\
\hline

{\bf II.}  &  $\mu=1$  & $ {}_{_sL}\mathcal{F}_n^{(m)}(y,x;\lambda;1,0)$  & 2-variable generalized Laguerre Apostol Euler polynomials\\

  &   $\nu=0$   &  $={}_{_sL}\mathfrak{E}_n^{(m)}(y,x;\lambda)$   & $\sum\limits_{k=0}^{n}{n\choose k}c^{n-k}d^{\nu+k} {}_{_sL}\mathfrak{E}_{n-k}^{(m)}(d^{s}y,dx;\lambda)\sum\limits_{l=0}^{k}{k\choose l}\mathcal{M}_{l}(c-1;\lambda){}_{_sL}\mathfrak{E}_{k-l}^{(m-1)}(c^{s}Y,cX;\lambda)$\\

     &    &    &$=\sum\limits_{k=0}^{n}{n\choose k}d^{n-k}c^{\nu+k} {}_{_sL}\mathfrak{E}_{n-k}^{(m)}(c^{s}y,cx;\lambda)\sum\limits_{l=0}^{k}{k\choose l}\mathcal{M}_{l}(d-1;\lambda){}_{_sL}\mathfrak{E}_{k-l}^{(m-1)}(d^{s}Y,dX;\lambda)$  \\
&&&\\
  &     &        & $\sum\limits_{k=0}^n {n \choose k} \sum\limits_{i=0}^{c-1}\sum\limits_{j=0}^{d-1} (-1)^{i+j}\lambda^{i+j} c^k d^{n-k} {}_{_sL}\mathfrak{E}_n^{(m)}\Big(d^{s}y,dx+\frac{d}{c}i;\lambda){}_{_sL}\mathfrak{E}_n^{(m)}\Big(c^{s}Y,cX+\frac{c}{d}j;\lambda\Big)$ \\

     &    &     &  $\sum\limits_{k=0}^n {n \choose k} \sum\limits_{i=0}^{d-1}\sum\limits_{j=0}^{c-1}(-1)^{i+j} \lambda^{i+j} d^k c^{n-k} {}_{_sL}\mathfrak{E}_n^{(m)}\Big(c^{s}y,cx+\frac{c}{d}i;\lambda\Big){}_{_sL}\mathfrak{E}_n^{(m)}\Big(d^{s}Y,dX+\frac{d}{c}j;\lambda\Big)$  \\
\hline
\end{tabular}}}\\

{\tiny{
\begin{tabular}{|l|l|l|l|l|}
\hline
&&&\\
{\bf III.}  &  $\mu=1$  & $ {}_{_sL}\mathcal{F}_n^{(m)}(y,x;\lambda;1,1)$  & 2-variable generalized Laguerre Apostol Genocchi polynomials\\

  &   $\nu=1$   &  $={}_{_sL}\mathcal{G}_n^{(m)}(y,x;\lambda)$   & $\sum\limits_{k=0}^{n}{n\choose k}c^{n-k}d^{\nu+k} {}_{_sL}\mathcal{G}_{n-k}^{(m)}(d^{s}y,dx;\lambda)\sum\limits_{l=0}^{k}{k\choose l}\mathcal{M}_{l}(c-1;\lambda){}_{_sL}\mathcal{G}_{k-l}^{(m-1)}(c^{s}Y,cX;\lambda)$ \\

     &    &    &$=\sum\limits_{k=0}^{n}{n\choose k}d^{n-k}c^{\nu+k} {}_{_sL}\mathcal{G}_{n-k}^{(m)}(c^{s}y,cx;\lambda)\sum\limits_{l=0}^{k}{k\choose l}\mathcal{M}_{l}(d-1;\lambda){}_{_sL}\mathcal{G}_{k-l}^{(m-1)}(d^{s}Y,dX;\lambda)$ \\
&&&\\
  &     &        & $\sum\limits_{k=0}^n {n \choose k} \sum\limits_{i=0}^{c-1}\sum\limits_{j=0}^{d-1} (-1)^{i+j}\lambda^{i+j} c^k d^{n-k} {}_{_sL}\mathcal{G}_n^{(m)}\Big(d^{s}y,dx+\frac{d}{c}i;\lambda){}_{_sL}\mathcal{G}_n^{(m)}\Big(c^{s}Y,cX+\frac{c}{d}j;\lambda\Big)$ \\

     &    &     &  $\sum\limits_{k=0}^n {n \choose k} \sum\limits_{i=0}^{d-1}\sum\limits_{j=0}^{c-1}(-1)^{i+j} \lambda^{i+j} d^k c^{n-k} {}_{_sL}\mathcal{G}_n^{(m)}\Big(c^{s}y,cx+\frac{c}{d}i;\lambda\Big){}_{_sL}\mathcal{G}_n^{(m)}\Big(d^{s}Y,dX+\frac{d}{c}j;\lambda\Big)$  \\
\hline

\end{tabular}}}\\
\vspace{.35cm}

\noindent
{\bf Note.}~We note that for $\lambda=1$, the symmetry derived above for the 2VGLABP $_{_sL}\mathfrak{B}_n^{(m)}(x,y;\lambda)$, 2VGLAEP $_{_sL}\mathfrak{E}_n^{(m)}(x,y;\lambda)$ and 2VGLAGP $_{_sL}\mathcal{G}_n^{(m)}(x,y;\lambda)$ give the corresponding identities for the 2-variable generalized Laguerre Bernoulli polynomials (2VGLBP) (of order $m$)  $_{_sL}B_n^{(m)}(x,y)$, 2-variable generalized Laguerre Euler polynomials (2VGLEP) (of order $m$) $_{_sL}E_n^{(m)}(x,y)$ and 2-variable generalized Laguerre Genocchi polynomials (2VGLGP) (of order $m$) $_{_sL}G_n^{(m)}(x,y)$, respectively. Again, for $m=1$, we get the corresponding identities for the 2-variable generalized Laguerre Bernoulli polynomials (2VGLBP) $_{_sL}B_n(x,y)$, 2-variable generalized Laguerre Euler polynomials (2VGLEP) $_{_sL}E_n(x,y)$ and 2-variable generalized Laguerre Genocchi polynomials (2VGLGP) $_{_sL}G_n(x,y)$.\\

Similarly, by taking $\phi(y,t)=\frac{1}{1-yt^r}$ (for which the 2VGP $p_n(x,y)$ reduce to the 2-variable truncated exponential polynomials (2VTE) $e_n^{(r)}(x,y))$ in the l.h.s. of generating function (1.15), we get the following generating function for the 2-variable truncated exponential Apostol type polynomials ${}_{e^{(r)}}\mathcal{F}_n^{(\alpha)}(x,y;\lambda;\mu,\nu)$ \cite{ATP}:

$$\Big(\frac{2^\mu~ t^\nu}{\lambda e^{t}+1}\Big)^m \Big(\frac{e^{xt}}{1-yt^r}\Big)=\sum\limits_{n=0}^{\infty} {}_{e^{(r)}}\mathcal{F}_n^{(m)}(x,y;\lambda;\mu,\nu)\frac{t^{n}}{n!}.\eqno(3.7)$$

Consequently, by taking $\phi(y,cdt)=\frac{1}{1-y(cdt)^r}$ in equation (2.2) and proceeding on the same lines as in Theorem 2.1, we get the following consequence of Theorem 2.1.\\

\noindent
{\bf Corollary 3.5.}~{\em For all integers $c>0,~d >0,~n \geq0, ~m \in \mathbb{N}, ~\lambda,\mu,\nu \in \mathbb{C}$, the following symmetry identity for the 2-variable truncated exponential Apostol type polynomials $_{e^{(r)}}\mathcal{F}_n^{(m)}(x,y;\lambda;\mu,\nu)$ holds true:}
\[\begin{split}
&\sum_{k=0}^{n}\sum_{l=0}^{k}{n\choose k}{k\choose l}c^{n-k}d^{\nu+k} {}_{{e}^{(r)}}\mathcal{F}_{n-k}^{(m)}(dx,d^{r}y;\lambda;\mu,\nu)S_{l}(c-1;-\lambda){}_{{e}^{(r)}}\mathcal{F}_{k-l}^{(m-1)}(cX,c^{r}Y;\lambda;\mu,\nu)\\
&=\sum_{k=0}^{n}\sum_{l=0}^{k}{n\choose k}{k\choose l}d^{n-k}c^{\nu+k} {}_{{e}^{(r)}}\mathcal{F}_{n-k}^{(m)}(cx,c^{r}y;\lambda;\mu,\nu)S_{l}(d-1;-\lambda){}_{{e}^{(r)}}\mathcal{F}_{k-l}^{(m-1)}
(dX,d^{r}Y;\lambda;\mu,\nu).\\
&\hspace{5.14in}(3.8)
\end{split}\]

\noindent
{\bf Corollary 3.6.}~{\em For each pair of positive integers $c,~d$ and for all integers $n \geq0, ~m \in \mathbb{N}, ~\lambda,\mu,\nu \in \mathbb{C}$, the following symmetry identity for the 2-variable truncated exponential Apostol type polynomials $_{e^{(r)}}\mathcal{F}_n^{(m)}(x,y;\lambda;\mu,\nu)$ holds true}:
\[\begin{split}
&\sum\limits_{k=0}^n {n \choose k} \sum\limits_{i=0}^{c-1}\sum\limits_{j=0}^{d-1} (-1)^{i+j}\lambda^{i+j} c^k d^{n-k} {}_{e^{(r)}}\mathcal{F}_n^{(m)}\Big(dx+\frac{d}{c}i,d^{r}y;\lambda;\mu,\nu\Big){}_{e^{(r)}}\mathcal{F}_n^{(m)}\Big(cX+\frac{c}{d}j,c^{r}Y;\lambda;\mu,\nu\Big)\\
&=\sum\limits_{k=0}^n {n \choose k} \sum\limits_{i=0}^{d-1}\sum\limits_{j=0}^{c-1}(-1)^{i+j} \lambda^{i+j} d^k c^{n-k} {}_{e^{(r)}}\mathcal{F}_n^{(m)}\Big(cx+\frac{c}{d}i,c^{r}y;\lambda;\mu,\nu\Big){}_{e^{(r)}}
\mathcal{F}_n^{(m)}\Big(dX+\frac{d}{c}j,d^{r}Y;\lambda;\mu,\nu\Big).\\
&\hspace{5.14in}(3.9)
\end{split}\]

By taking suitable values of the parameters in equation (3.2) and (3.3) and in view of relations (1.12)-(1.14), we can find the symmetric identities for other mixed special polynomials related to $_{e^{(r)}}\mathcal{F}_n^{(m)}(x,y;\lambda;\mu,\nu)$. We present the symmetric identities for these polynomials in the following table:\\
\vspace{.35cm}

\noindent
{\bf Table~3.3.~Symmetry identities for the Special cases of the 2VTEATP $_{e^{(r)}}\mathcal{F}_n^{(m)}(x,y;\lambda;\mu,\nu)$}\\
\\
{\tiny{
\begin{tabular}{|l|l|l|l|l|}
\hline
&&&\\
{\bf S.} & {\bf Values of}  & {\bf Relation between}  &  {\bf Name and Symmetry identities  of the resultant special polynomials}\\

 {\bf No.} &  {\bf the para-}  &  {\bf the 2VTEATP}   & \\
   & {\bf meters}  &  $_{e^{(r)}}\mathcal{F}_n^{(m)}(x,y;\lambda;\mu,\nu)$      &  \\
  &    &  {\bf and its special case}    &   \\
\hline
&&&\\
{\bf I.}  &  $\lambda \rightarrow -\lambda$  & $(-1)^m $  & 2-variable truncated exponential Apostol Bernoulli polynomials \\

  &  $\mu=0$   & ${}_{e^{(r)}}\mathcal{F}_n^{(m)}(x,y;-\lambda;0,1)$   &  $\sum\limits_{k=0}^{n}{n\choose k}c^{n-k}d^{\nu+k} {}_{e^{(r)}}\mathfrak{B}_{n-k}^{(m)}(dx,d^{r}y;\lambda)\sum\limits_{l=0}^{k}{k\choose l}\mathcal{S}_{l}(c-1;\lambda){}_{e^{(r)}}\mathfrak{B}_{k-l}^{(m-1)}(cX,c^{r}Y;\lambda)$ \\

     &   $\nu=1$  &  $={}_{e^{(r)}}\mathfrak{B}_n^{(m)}(x,y;\lambda)$  & $=\sum\limits_{k=0}^{n}{n\choose k}d^{n-k}c^{\nu+k} {}_{e^{(r)}}\mathfrak{B}_{n-k}^{(m)}(cx,c^{r}y;\lambda)\sum\limits_{l=0}^{k}{k\choose l}\mathcal{S}_{l}(d-1;\lambda){}_{e^{(r)}}\mathfrak{B}_{k-l}^{(m-1)}(dX,d^{r}Y;\lambda)$\\
&&&\\
  &     &         & $\sum\limits_{k=0}^n {n \choose k} \sum\limits_{i=0}^{c-1}\sum\limits_{j=0}^{d-1} (-1)^{i+j}(-\lambda)^{i+j} c^k d^{n-k} {}_{e^{(r)}}\mathfrak{B}_n^{(m)}\Big(dx+\frac{d}{c}i,d^{r}y;\lambda){}_{e^{(r)}}\mathfrak{B}_n^{(m)}\Big(cX+\frac{c}{d}j,c^{r}Y;\lambda\Big)$ \\

     &    &     &  $\sum\limits_{k=0}^n {n \choose k} \sum\limits_{i=0}^{d-1}\sum\limits_{j=0}^{c-1}(-1)^{i+j} (-\lambda)^{i+j} d^k c^{n-k} {}_{e^{(r)}}\mathfrak{B}_n^{(m)}\Big(cx+\frac{c}{d}i,c^{r}y;\lambda\Big){}_{e^{(r)}}\mathfrak{B}_n^{(m)}\Big(dX+\frac{d}{c}j,d^{r}Y;\lambda\Big)$  \\
\hline

{\bf II.}  &  $\mu=1$  & $ {}_{e^{(r)}}\mathcal{F}_n^{(m)}(x,y;\lambda;1,0)$  & 2-variable truncated exponential Apostol Euler polynomials\\

  &   $\nu=0$   &  $={}_{e^{(r)}}\mathfrak{E}_n^{(m)}(x,y;\lambda)$   & $\sum\limits_{k=0}^{n}{n\choose k}c^{n-k}d^{\nu+k} {}_{e^{(r)}}\mathfrak{E}_{n-k}^{(m)}(dx,d^{r}y;\lambda)\sum\limits_{l=0}^{k}{k\choose l}\mathcal{M}_{l}(c-1;\lambda){}_{e^{(r)}}\mathfrak{E}_{k-l}^{(m-1)}(cX,c^{r}Y;\lambda)$\\

     &    &    &$=\sum\limits_{k=0}^{n}{n\choose k}d^{n-k}c^{\nu+k} {}_{e^{(r)}}\mathfrak{E}_{n-k}^{(m)}(cx,c^{r}y;\lambda)\sum\limits_{l=0}^{k}{k\choose l}\mathcal{M}_{l}(d-1;\lambda){}_{e^{(r)}}\mathfrak{E}_{k-l}^{(m-1)}(dX,d^{r}Y;\lambda)$  \\
&&&\\
  &     &        & $\sum\limits_{k=0}^n {n \choose k} \sum\limits_{i=0}^{c-1}\sum\limits_{j=0}^{d-1} (-1)^{i+j}\lambda^{i+j} c^k d^{n-k} {}_{e^{(r)}}\mathfrak{E}_n^{(m)}\Big(dx+\frac{d}{c}i,d^{r}y;\lambda){}_{e^{(r)}}\mathfrak{E}_n^{(m)}\Big(cX+\frac{c}{d}j,c^{r}Y;\lambda\Big)$ \\

     &    &     &  $\sum\limits_{k=0}^n {n \choose k} \sum\limits_{i=0}^{d-1}\sum\limits_{j=0}^{c-1}(-1)^{i+j} \lambda^{i+j} d^k c^{n-k} {}_{e^{(r)}}\mathfrak{E}_n^{(m)}\Big(cx+\frac{c}{d}i,c^{r}y;\lambda\Big){}_{e^{(r)}}\mathfrak{E}_n^{(m)}\Big(dX+\frac{d}{c}j,d^{r}Y;\lambda\Big)$  \\
\hline

{\bf III.}  &  $\mu=1$  & $ {}_{e^{(r)}}\mathcal{F}_n^{(m)}(x,y;\lambda;1,1)$  & 2-variable truncated exponential Apostol Genocchi polynomials\\

  &   $\nu=1$   &  $={}_{e^{(r)}}\mathcal{G}_n^{(m)}(x,y;\lambda)$   & $\sum\limits_{k=0}^{n}{n\choose k}c^{n-k}d^{\nu+k} {}_{e^{(r)}}\mathcal{G}_{n-k}^{(m)}(dx,d^{r}y;\lambda)\sum\limits_{l=0}^{k}{k\choose l}\mathcal{M}_{l}(c-1;\lambda){}_{e^{(r)}}\mathcal{G}_{k-l}^{(m-1)}(cX,c^{r}Y;\lambda)$ \\

     &    &    &$=\sum\limits_{k=0}^{n}{n\choose k}d^{n-k}c^{\nu+k} {}_{e^{(r)}}\mathcal{G}_{n-k}^{(m)}(cx,c^{r}y;\lambda)\sum\limits_{l=0}^{k}{k\choose l}\mathcal{M}_{l}(d-1;\lambda){}_{e^{(r)}}\mathcal{G}_{k-l}^{(m-1)}(dX,d^{r}Y;\lambda)$ \\
&&&\\
  &     &        & $\sum\limits_{k=0}^n {n \choose k} \sum\limits_{i=0}^{c-1}\sum\limits_{j=0}^{d-1} (-1)^{i+j}\lambda^{i+j} c^k d^{n-k} {}_{e^{(r)}}\mathcal{G}_n^{(m)}\Big(dx+\frac{d}{c}i,d^{r}y;\lambda){}_{e^{(r)}}\mathcal{G}_n^{(m)}\Big(cX+\frac{c}{d}j,c^{r}Y;\lambda\Big)$ \\

     &    &     &  $\sum\limits_{k=0}^n {n \choose k} \sum\limits_{i=0}^{d-1}\sum\limits_{j=0}^{c-1}(-1)^{i+j} \lambda^{i+j} d^k c^{n-k} {}_{e^{(r)}}\mathcal{G}_n^{(m)}\Big(cx+\frac{c}{d}i,c^{r}y;\lambda\Big){}_{e^{(r)}}\mathcal{G}_n^{(m)}\Big(dX+\frac{d}{c}j,d^{r}Y;\lambda\Big)$  \\
\hline

\end{tabular}}}\\
\vspace{.35cm}

\noindent
{\bf Note.}~We note that for $\lambda=1$, the symmetry derived above for the 2VTEABP $_{e^{(r)}}\mathfrak{B}_n^{(m)}(x,y;\lambda)$, 2VTEAEP $_{e^{(r)}}\mathfrak{E}_n^{(m)}(x,y;\lambda)$ and 2VTEAGP $_{e^{(r)}}\mathcal{G}_n^{m)}(x,y;\lambda)$ give the corresponding identities for the 2-variable truncated exponential Bernoulli polynomials (2VTEBP) (of order $\alpha$)  $_{e^{(r)}}B_n^{(\alpha)}(x,y)$, 2-variable truncated exponential Euler polynomials (2VTEEP) (of order $m$) $_{e^{(r)}}E_n^{(m)}(x,y)$ and 2-variable truncated exponential Genocchi polynomials (2VTEGP) (of order $m$) $_{e^{(r)}}G_n^{(m)}(x,y)$, respectively. Again, for $m=1$, we get the corresponding identities for the 2-variable truncated exponential Bernoulli polynomials (2VTEBP) $_{e^{(r)}}B_n(x,y)$, 2-variable truncated exponential Euler polynomials (2VTEEP) $_{e^{(r)}}E_n(x,y)$ and 2-variable truncated exponential Genocchi polynomials (2VTEGP) $_{e^{(r)}}G_n(x,y)$.\\

\end{document}